\documentclass[12pt]{article}
\usepackage{amssymb,amsmath}
\usepackage{mathrsfs}
\usepackage{graphicx}

\newtheorem{Thm}{Theorem}
\newtheorem{Lm}{Lemma}

\newcommand{\sgn}{\mathrm{sgn}}

\makeatletter
\renewcommand{\@cite}[2]{{#1\if@tempswa , #2\fi}}
\renewcommand{\@biblabel}[1]{\hfill}
\makeatother

\begin{document}

\begin{center}
\large{\bf Catalytic Branching Random Walk\\
with Semi-exponential Increments}
\end{center}
\vskip0,5cm
\begin{center}
Ekaterina Vl. Bulinskaya\footnote{ \emph{Email address:} {\tt
bulinskaya@yandex.ru}}$^,$\footnote{The work is supported by Russian Science Foundation under grant 17-11-01173 and is fulfilled
at Novosibirsk State University. The author is Associate Professor of the Lomonosov Moscow State University.}
\vskip0,2cm \emph{Novosibirsk State University}
\end{center}

\begin{abstract}

A catalytic branching random walk on a multidimensional lattice, with arbitrary finite number of catalysts, is studied in supercritical regime. The dynamics of spatial spread of the particles population is examined, upon normalization. The components of the vector random walk jump are assumed independent (or close to independent) and have semi-exponential distributions with, possibly, different parameters. A  limit theorem on the almost sure normalized positions of the particles at the population ``front'' is established. Contrary to the case of the random walk increments with ``light'' distribution tails, studied by Carmona and Hu (2014) in one-dimensional setting and Bulinskaya (2018) in multidimensional setting, the normalizing factor has a power rate and grows faster than linear in time function. The limiting shape of the front in the case of semi-exponential tails is non-convex in contrast to a convex one in the case of light tails.

\vskip0,5cm {\it Keywords and phrases}: catalytic branching random walk, supercritical regime,
heavy tails, semi-exponential distribution tails, spread of population, population front.

\vskip0,5cm 2010 {\it AMS classification}: 60J80, 60F15.

\end{abstract}

\section{Introduction}

Stochastic models describing evolution of certain arrays of particles or  biological species (genes, bacteria, infected individuals) are of theoretical interest and also important for various applications, \cite{New_14}, \cite{Bansaye_15}, \cite{Mel_16}, \cite{Pardoux_16}, \cite{KV_17}. Investigation of time-dependence of population size in probabilistic terms ascends to introduction of the classical Galton-Watson branching process in 1874, see \cite{Jagers_11} for historical overview. A rigorous analysis of  the spread of a new dominant gene can be traced back to the well-known paper by \cite{KPP_37}. Later a plenty of models involving reproduction, death, and movement of particles have appeared. The probabilistic ones include the so-called branching random walks (BRW) and are linked with other models in mathematics, physics, and informatics. Most of BRW models are space-homogeneous as in \cite{Lifshits_14}, \cite{Shi_LNM_15}, \cite{Mallein_16}, and references therein. We consider a non-homogeneous process, called \emph{catalytic branching random walk} (CBRW), defined for any time $t\geq 0$ on a multidimensional lattice $\mathbb{Z}^d$, $d\in\mathbb{N}$. Particles give offspring or die at specified locations called \emph{catalysts}. The catalysts take fixed positions on $\mathbb{Z}^d$. Their number is arbitrary finite. Outside the catalysts the particles move randomly until hitting a catalyst.

So CBRW models incorporate two stochastic mechanisms: the  particles randomly move in space and, moreover, each of them could give
a random number of offsprings in the presence of catalysts only. These issues are discussed further in detail for the populations
initiated by a single specimen. Our main goal is to examine the time-evolution of the moving front, separating, in a sense,
the once populated area and its empty environment. Whenever the locations of particles are normalized in appropriate manner by
continuous in time functions, the rescaled front is a surface in $\mathbb{R}^d$. We study the spread of population for $t\to\infty$
when the distributions of particles jumps have heavy tails. The exact limiting front surface in $\mathbb{R}^d$, called the limiting shape of the front, is identified with probability one. The proofs combine analysis of non-linear integral equations, multidimensional renewal theorems, the Laplace transform, large deviation theory for heavy-tailed distributions, the coupling method, and other probabilistic-analytic techniques.

Likewise the classical branching processes (\cite{Sev_74}), CBRW is classified as supercritical, critical, and subcritical depending
on the characteristics of the reproduction and of the random walk. In the most general framework the classification is given
by \cite{B_TPA_15} by means of the Perron eigenvalue of a certain matrix. The particles population survives globally and locally
with positive probability if and only if CBRW is supercritical, \cite{B_Doklady_15}. Exponential growth of the total and local particles numbers occurs only in supercritical CBRW as established in \cite{B_TPA_15} and \cite{B_Doklady_15}. This is the reason to consider the rate of population spread just for supercritical CBRW.

\cite{Carmona_Hu_14} study for CBRW on $\mathbb{Z}$ the strong (that is almost sure) limit behavior of the maximum, being the location of the right-most particle on $\mathbb{Z}$, under assumption of light distribution tails of the random walk. \cite{B_SPA_18} extends analysis to CBRW on the lattice of arbitrary dimension. Until now the investigation of the maximum of CBRW with heavy distribution tails of the random walk has been tackled only in \cite{B_Arxiv_18}, where the distribution tail is a regularly varying function. It follows from \cite{Carmona_Hu_14} that the unrescaled population front propagates linearly in time in case of light distribution tails whereas according to \cite{B_Arxiv_18} it spreads exponentially fast in case of regularly varying tails.

We focus here on a novel assumption for CBRW models that the distribution tails of the random walk are semi-exponential. Such distribution belongs to a fundamental class of heavy-tailed distributions including the Weibull one (\cite{Borovkov_Borovkov_08}, Ch.~5). It follows that the non-normalized population front on $\mathbb{Z}^d$ for CBRW, with independent components of the random walk jump distributed semi-exponentially, propagates in a power way and faster than a linear in time function.

An equation of the form $H({\bf z})=\nu$, ${\bf z}\in \mathbb{R}^d$, is obtained, describing the surface of the normalized front limit (limiting shape of the front) for the model under consideration, where $\nu$ is a Maltusian parameter, $H$ is an explicit function depending only on ${\bf z}$ and parameters of the semi-exponential distributions of the jumps components. The growth rate of the normalizing factors also depends exceptionally on these parameters.

Surprisingly, in the case of semi-exponential tails the limiting front shape is a surface of non-convex set which contrasts with a convex one in the case of light tails. Our results for CBRW on $\mathbb{Z}$ agree with those for \emph{homogeneous} branching random walk on the real line proved in \cite{Gantert_00}. One can refer to \cite{Maillard_16}, \cite{GMV_2017},
and others on progress in studying the spread of the population for homogeneous branching random walk with regularly varying jump distribution tails ``heavier'' than the semi-exponential ones.
We found only the paper \cite{Gantert_00} devoted to homogeneous branching random walk on the real
line with semi-exponential increments. There are no investigations of such model in a
multidimensional case. Possibly, it is explained by the absence of quite general results
on large deviations of vectors with semi-exponential distribution. In its turn, according
to \cite{Borovkov_Borovkov_08}, p.~400, the latter is a hard problem.
It seems that our work opens the study of multidimensional BRW models with semi-exponential jumps distributions.

Observe that CBRW captures the nature of the intermittency phenomenon for random walk in random potential (\cite{Konig_16}, p.10): the main contribution to the total population size is due to  few small remote islands, called intermittent islands. In CBRW the counterparts of the intermittent islands are the catalysts points. Investigations in mathematical theory of intermittency and applications to magnetic and temperature fields of turbulent flows, chemical kinetics, hydrodynamics, and biological models are performed in \cite{Gartner_Molchanov_90}, \cite{Zeldovich_14}, \cite{Ortgiese_Roberts_16}, \cite{Chernousova_Molchanov_18}.

We refer to the works by \cite{AB_00}, \cite{Molchanov_Yarovaya_12}, \cite{HTV_12}, \cite{Topchii_Vatutin_13}, \cite{DR_13}, and \cite{Platonova_Ryadovkin_17}, and also to references therein, for analysis of other aspects of CBRW or its modifications. Most of them are devoted to long-time behavior of total and local particles numbers. The exception is \cite{Molchanov_Yarovaya_12}, where the population front of symmetric CBRW with binary splitting and light-tailed increments was defined and studied from the viewpoint of moments boundedness of local particles numbers. Remarkably, the notions of the propagation front in \cite{Molchanov_Yarovaya_12} and \cite{B_SPA_18} are different but lead to the same growth rate. However, our approach seems more powerful due to the strong convergence results under milder restrictions on CBRW.

In Section~\ref{s:main_results_semiexp} we introduce necessary notation and formulate the main results in multidimensional setting.
Theorems~\ref{T:main_result_semiexp_lattice} and \ref{T:one_point_semiexp} characterize the front propagation almost surely.
In Section~\ref{s:proof_2_semiexp} we provide the proofs of these theorems. To simplify exposition we establish 6 lemmas and  partition the proofs into several steps. At first we consider the case of a single catalyst and then extend the obtained results to the case of an arbitrary finite number of catalysts. Illustrative examples are gathered in Section~\ref{s:examples_semiexp}.
Section Conclusion completes the exposition.

\section{Model Description and Main Results}
\label{s:main_results_semiexp}

All random elements are defined on a complete probability space $(\Omega,\mathcal{F},{\sf P})$, letter $\omega$ stands
for a point of $\Omega$. The index ${\bf z}$ in expressions of the form ${\sf E}_{\bf z}$ and ${\sf P}_{\bf z}$ marks the starting point of either CBRW or the random walk ${\bf S}$, depending on the context. Bold font of ${\bf z}$ emphasizes that ${\bf z}$ is a multidimensional vector whereas the symbol $z$ means that $z$ is a real number.

Recall the description of CBRW on $\mathbb{Z}^d$, $d\in\mathbb{N}$ (in our setting given in \cite{B_SPA_18}). At the initial time ${t=0}$ there is a single particle that moves on $\mathbb{Z}^d$ according to a continuous-time Markov chain ${\bf S}=\{{\bf S}(t),t\geq0\}$ generated by the infinitesimal matrix ${Q=(q({\bf x},{\bf y}))_{{\bf x},{\bf y}\in\mathbb{Z}^d}}$. Assume that ${\bf S}$ is irreducible and space-homogeneous, with the conservative matrix $Q$, that is,
$Q$ has finite elements and
\begin{equation}\label{condition1}
q({\bf x},{\bf y})=q({\bf x}-{\bf y},{\bf 0})
\quad\mbox{and}\quad
\sum\limits_{{\bf y}\in\mathbb{Z}^d}{q({\bf x},{\bf y})}=0,
\end{equation}
where $q({\bf x},{\bf y})\geq0$, for ${\bf x}\neq{\bf y}$, and $q:=-q({\bf x},{\bf x})\in(0,\infty)$, for any ${\bf x},{\bf y}\in\mathbb{Z}^d$. Stress that, contrary to Platonova and Ryadovkin (2017) and Yarovaya (2017), we do not restrict ourselves to the case of symmetric generator $Q$.

When this particle hits a finite set of catalysts $W=\{{\bf w}_1,\ldots,{\bf w}_N\}\subset\mathbb{Z}^d$, say at the point ${\bf w}_k$, it spends there random time, distributed exponentially with parameter $\beta_k>0$. The particle either branches there with probability $\alpha_k$ or leaves the point ${\bf w}_k$ with probability $1-\alpha_k$ ($0\leq\alpha_k<1$). If the particle branches, it produces a random non-negative integer number $\xi_k$ of offsprings, located at the same point ${\bf w}_k$, and dies instantly. Whenever the particle leaves ${\bf w}_k$, it jumps to the point ${\bf y}\neq{\bf w}_k$ with probability $-(1-\alpha_k)q({\bf w}_k,{\bf y})q({\bf w}_k,{\bf w}_k)^{-1}$ and resumes its motion governed by the Markov chain ${\bf S}$. All the newly born particles are supposed to behave as independent copies of their parent.

Denote by $f_k(s):={\sf E}{s^{\xi_k}}$, $s\in[0,1]$, the probability generating function of $\xi_k$, $k=1,\ldots,N$. Employ the standard assumption of existence of a finite derivative $f_k'(1)$, that is the finiteness of $m_k:={\sf E}{\xi_k}$, for any $k=1,\ldots,N$.

We forget for a while about the catalysts and consider only the motion of a particle on $\mathbb{Z}^d$ according to the Markov chain ${\bf S}$ with the generator $Q$ and the starting point ${\bf x}$. The conditions imposed on the elements $q({\bf x},{\bf y})$, ${\bf x},{\bf y}\in\mathbb{Z}^d$, allow us to use an explicit construction of the random walk on $\mathbb{Z}^d$ with generator $Q$ by Theorem~1.2 in \cite{Bremaud_99}, Ch.~9, Sec.~1. Whence ${\bf S}$ is a regular jump process with right continuous trajectories and, for transition times of the process, $\tau^{(0)}:=0$ and $\tau^{(n)}:=\inf\left\{t\geq\tau^{(n-1)}:{\bf S}(t)\neq {\bf S}(\tau^{(n-1)})\right\}$, $n\ge1$, the following property is valid. The random variables $\left\{\tau^{(n+1)}-\tau^{(n)}\right\}_{n=0}^{\infty}$ are independent and each of them has exponential distribution with parameter $q$. Denote by $\Pi=\{\Pi(t),t\geq0\}$ the Poisson process constructed as the renewal process with the interarrival times $\tau^{(n+1)}-\tau^{(n)}$, $n\in\mathbb{Z}_+$, (\cite{Feller_71}, Ch.~1, Sec.~4), that is, $\Pi$ is a Poisson process with constant intensity $q$. Let ${\bf Y}^i=\left(Y^i_1,\ldots,Y^i_d\right)$ be the value of the $i$th jump of the random walk ${\bf S}$ ($i=1,2,\ldots$). In view of Theorem~1.2 in \cite{Bremaud_99}, Ch.~9, Sec.~1, the random vectors ${\bf Y}^1,{\bf Y}^2,\ldots$ are i.i.d., have distribution ${\sf P}({\bf Y}^1={\bf y})=q({\bf 0},{\bf y})/q$, ${\bf y}\in\mathbb{Z}^d$, ${\bf y}\neq{\bf 0}$, and do not depend on the sequence $\{\tau^{(n+1)}-\tau^{(n)}\}_{n=0}^{\infty}$. We can write (for a version of ${\bf S}$)
\begin{equation}\label{S(t)=representation}
{\bf S}(t)={\bf x}+\sum_{i=1}^{\Pi(t)}{\bf Y}^i,
\end{equation}
where  ${\bf x}$ is the initial state of the Markov chain ${\bf S}$ and  $\sum_{i\in\varnothing}{\bf Y}^i:={\bf 0}$. Equality~(\ref{S(t)=representation}) implies that ${\bf S}$ is a process with independent increments. In what follows we consider the version of the process ${\bf S}$ constructed in such a way. It is called a compound Poisson process.

We employ various stopping times (with respect to the natural filtration of ${\bf S}$).
For ${\bf x}\in \mathbb{Z}^d$, set
$$\tau_{\bf x}:=\mathbb{I}({\bf S}(0)={\bf x})\inf\{t\geq0:{\bf S}(t)\neq{\bf x}\},$$
that is, the stopping time $\tau_{\bf x}$ is the time of the first exit from the starting point ${\bf x}$ of the random walk. As usual,
$\mathbb{I}(A)$ stands for the indicator of a set $A\in\mathcal{F}$. Clearly, ${\sf P}_{\bf x}(\tau_{\bf x}\leq t)=1-e^{-qt}=:G_0(t)$, ${\bf x}\in\mathbb{Z}^d$, $t\geq0$. Let
$$_T\overline{\tau}_{{\bf x},{\bf y}}:=\mathbb{I}({\bf S}(0)={\bf x})\inf\{t\geq0:{\bf S}(t+\tau_{\bf x})={\bf y},{\bf S}(u)\notin T,\tau_{\bf x}\leq u< t+\tau_{\bf x}\}$$
be the time elapsed from the exit moment of this Markov chain (in other terms, particle) from the starting point
${\bf x}$ till the moment of the first hitting point ${\bf y}$, whenever the particle trajectory does not pass the set $T\subset\mathbb{Z}^d$. If there is no such finite $t$, we put ${_T\overline{\tau}_{{\bf x},{\bf y}}=\infty}$. An extended random variable $_T\overline{\tau}_{{\bf x},{\bf y}}$ is called \emph{hitting time} of state ${\bf y}$ \emph{under taboo} on set $T$ after exit out of starting state ${\bf x}$ (\cite{Chung_60}, Ch.~2, Sec.~11, and \cite{B_SPL_14}). Denote by $_T\overline{F}_{{\bf x},{\bf y}}(t)$, $t\geq0$, the improper c.d.f. of this extended random variable and let $_T\overline{F}_{{\bf x},{\bf y}}(\infty):=\lim_{t\to\infty}{_T\overline{F}_{{\bf x},{\bf y}}(t)}$. If the taboo set $T$ is empty, expressions $_{\varnothing}\overline{\tau}_{{\bf x},{\bf y}}$ and $_{\varnothing}\overline{F}_{{\bf x},{\bf y}}$ are shortened as $\overline{\tau}_{{\bf x},{\bf y}}$ and $\overline{F}_{{\bf x},{\bf y}}$. Mainly we will be interested in the situation when $T=W_k$, where $W_k:=W\setminus\{{\bf w}_k\}$, $k=1,\ldots,N$.

We  also use some auxiliary functions. Here and further let
$$F^{\ast}(\lambda):=\int\nolimits_{0-}^{\infty}{ e^{-\lambda t}\,d{F(t)}},\quad\lambda\geq0,$$
denote the Laplace transform of a c.d.f. $F(t)$, $t\geq0$, with support on non-negative semi-axis. For $j,k=1,\ldots,N$, ${\bf x},{\bf y}\in\mathbb{Z}^d$, and $t\geq0$, set
\begin{equation}\label{G-i,G-ij=}
G_j(t):=1-e^{-\beta_j t},\quad G_{j,k}(t):=G_j\ast{_{W_k}\overline{F}_{{\bf w}_j,{\bf w}_k}(t)},\quad{_T F_{{\bf x},{\bf y}}(t)}:=G_0\ast{_T\overline{F}_{{\bf x},{\bf y}}(t)},
\end{equation}
where $\ast$ between c.d.f. stands for their convolution. By definition the function ${_T F_{{\bf x},{\bf y}}(\cdot)}$ is a c.d.f. of the variable $_T\tau_{{\bf x},{\bf y}}:=\tau_{\bf x}+{_T\overline{\tau}_{{\bf x},{\bf y}}}$ called \emph{hitting time} of state ${\bf y}$ \emph{under taboo} on set $T$ when the starting state is ${\bf x}$.

As in \cite{B_TPA_15}, consider a matrix function $D(\lambda)=\left(d_{i,j}(\lambda)\right)_{i,j=1}^N$, $\lambda\geq0$, taking values in the set of irreducible matrices of size $N\times N$, with elements defined by way of
$$d_{i,j}(\lambda)=\delta_{i,j}\alpha_im_iG^{\ast}_i(\lambda)+(1-\alpha_i)G^{\ast}_{i,j}(\lambda),$$ where $\delta_{i,j}$ is the Kronecker delta. According to Definition~$1$ in \cite{B_TPA_15} CBRW is called {\it supercritical} if the Perron root (that is, positive eigenvalue being the spectral radius) $\rho(D(0))$ of the matrix $D(0)$ is greater than $1$. In view of monotonicity of
all elements of the matrix function $D(\cdot)$ there exists the solution $\nu>0$ of the equation $\rho(D(\lambda))=1$. On account of Theorem $1$ in \cite{B_TPA_15} just this positive number $\nu$ specifies the rate of exponential growth of the mean total and local particles numbers (in the literature devoted to population dynamics and classical branching processes one traditionally speaks of {\it Malthusian parameter}). In the sequel we consider the supercritical CBRW on $\mathbb{Z}^d$.

Let $N(t)\subset\mathbb{Z}^d$ be the (random) set of particles existing in CBRW at time ${t\geq0}$. For a particle $v\in N(t)$, denote by ${\bf X}^v(t)=\left(X^v_1(t),\ldots,X^v_d(t)\right)$ its position at time $t$. Consider the set
$$
\mathcal{I}:=\left\{\omega:\limsup_{t\to\infty}\{v\in N(t):{\bf X}^v(t)\in W\}\neq\varnothing\right\}\in\mathcal{F}.
$$
To avoid operations with a continuum number of sets $\left\{A_t\right\}_{t\geq0}$, we just put $\limsup_{t\to\infty}A_t:=\cap_{m=1}^{\infty}\cup_{k=1}^{\infty}\cap_{n=k}^{\infty}A_{n/2^m}$, that is, we deal with the binary rational values of the parameter $t$ only instead of its all non-negative values. For each $\omega\in\mathcal{I}$, there is an increasing to infinity sequence of \emph{binary rational} values $t_l(\omega)$, $l\in\mathbb{N}$, such that at each time $t_l(\omega)$ there are particles at $W$. The event consisting of $\omega$ for which there exists a similar sequence of \emph{any} (not only binary rational) values $t_l(\omega)$ is of the same probability ${\sf P}(\mathcal{I})$, and we may call $\mathcal{I}$ {\it the event of infinite number of visits of catalysts}. The behavior of CBRW on the set complement $\overline{\mathcal{I}}$ is a.s. trivial. Indeed, for $t\geq t_0(\omega)$ large enough either CBRW dies out or CBRW constitutes the system of some random walks (without branching) starting respectively from ${\bf X}^v(\omega,t_0)$, $v\in N(t_0)$, at time $t_0$. The supercritical regime of CBRW guarantees that ${\sf P}(\mathcal{I})>0$ (Theorem~4 of \cite{B_Doklady_15}).

Assumptions made are of the same type as in the previous papers \cite{Carmona_Hu_14} (with discrete  time on $\mathbb{Z}$), \cite{B_SPA_18}, and \cite{B_Arxiv_18} devoted to the study of the population spread in CBRW.
The following hypothesis corresponds to the new case under consideration.
Let the components of the random walk jumps be semi-exponentially distributed,
that is, for any $i=1,\ldots,d$ and $y\in\mathbb{Z}_+$, one has by  \cite{Borovkov_Borovkov_08}, p.~29,
\begin{eqnarray}
{\sf P}(Y^1_i>y)&=&L^{(1,+)}_i(y)\exp\left\{-y^{\gamma^{+}_i}L^{(2,+)}_i(y)\right\}:=R^{+}_i(y),\label{assumption:tails_right_semiexp}\\
{\sf P}(Y^1_i<-y)&=&L^{(1,-)}_i(y)\exp\left\{-y^{\gamma^{-}_i}L^{(2,-)}_i(y)\right\}:=R^{-}_i(y).\label{assumption:tails_left_semiexp}
\end{eqnarray}
Symbol ``$+$'' marks the right distribution tail whereas symbol ``$-$'' refers to the left one. For each $i=1,\ldots,d$ and $\kappa\in\{+,-\}$, functions $L^{(1,\kappa)}_i(y)$ and $L^{(2,\kappa)}_i(y)$, ${y\in\mathbb{Z}_+}$, are slowly varying, while parameters $\gamma^{\kappa}_i$ are taken from $(0,1)$. Recall that
$$
{\sf P}\left(Y^1_i>y\right)={q^{-1}\sum_{{\bf x}:\,x_i>y}{q({\bf 0},{\bf x})}},
$$
where ${\bf x}=\left(x_1,\ldots,x_d\right)\in\mathbb{Z}^d$.

It follows from (\ref{assumption:tails_right_semiexp}) and (\ref{assumption:tails_left_semiexp}) that $-\ln R^{\kappa}_i(y)$, $y\in\mathbb{Z}_+$, is a regularly varying function of index $\gamma^{\kappa}_i$. In accordance with \cite{Seneta_76}, Ch.~1, Sec.~5, property $5^\circ$, there exists an asymptotically uniquely determined inverse function $R^{-1,\,\kappa}_i(s)$, $s\geq0$, in the sense that $-\ln{R^{\,\kappa}_i\left(R^{-1,\,\kappa}_i(y)\right)}\sim y$, $R^{-1,\,\kappa}_i\left(-\ln{R^{\,\kappa}_i(y)}\right)\sim y$, as $y\to\infty$, $y\in\mathbb{Z}_+$, and
$$
R^{-1,\kappa}_i(s)=s^{1/\gamma^{\kappa}_i}L^{(3,\kappa)}_i(s),
$$
where $L^{(3,\kappa)}_i(s)$, $s\geq0$, is a slowly varying function at infinity.

Functions $R^{-1,\kappa}_i(\cdot)$, $i=1,\ldots,d$, play an important role in normalization of
coordinates of particles of $N(t)$. Emphasize that we have to use notation involving
$\kappa$ since the normalization of a particle position, in general, depends on
the orthant (one among $2^d$) of $\mathbb{R}^d$ containing this particle.

Unlike the random walks with either light or regularly varying distribution tails, a diversity of large deviations zones is
inherent in case of walks with semi-exponential increments. They are \emph{Cr\'{a}mer deviation} zone, \emph{intermediate} zone, and \emph{maximum jump approximation} zone, \cite{Borovkov_Borovkov_08}, p.~238. We deal with two latter ones. Writing in a compact form, assume that, for each fixed ${\bf x}=\left(x_1,\ldots,x_d\right)\in\mathbb{R}^d$, ${\bf x}\neq{\bf 0}$, one has
\begin{eqnarray}\label{assumptions:large_deviations_tails_right}
& &{\sf P}_{\bf 0}\left(\sgn({\bf x}){\bf S}(u)/{\bf R}^{-1,\kappa({\bf x})}(t)\in\left[\left|{\bf x}\right|,+\infty\right)\right)\nonumber\\
&=&h(u)\left(1+\delta(u,t)\right)\prod_{i=1}^d\left({\sf P}\left(\sgn(x_i)Y_i\geq|x_i| R^{-1,\kappa(x_i)}_i(t)\right)\right)^{\left(1+\varepsilon_i(t)\right)},
\end{eqnarray}
where $h(u)$, $u\geq0$, is a positive non-decreasing function such that $h(u)\sim c u^d$, $u\to\infty$, for
a constant $c>0$, the functions $\varepsilon_i(t)=\varepsilon_i(t,{\bf x})\to0$, as $t\to\infty$, and $\delta(u,t)=\delta(u,t,{\bf x})\to0$, as $u,t\to\infty$, $u\leq t$, $i=1,\ldots,d$. In relation (\ref{assumptions:large_deviations_tails_right}) the notation $\sgn({\bf x}){\bf S}(u)/{\bf R}^{-1,\kappa({\bf x})}(t)$ means the vector in $\mathbb{R}^d$ with $i$th component $\sgn(x_i)S_i(t)/R^{-1,\kappa(x_i)}_i(t)$, $i=1,\ldots,d$, and $[|{\bf x}|,+\infty):=[|x_1|,+\infty)\times\ldots\times[|x_d|,+\infty)$. Here $\kappa(x_i)=$``$+$'' if $x_i\geq0$ and $\kappa(x_i)=$``$-$'' if $x_i<0$. As a precaution, we put $\sgn(x_i)S_i(t)/R^{-1,\kappa(x_i)}_i(t):=0$ whenever $R^{-1,\kappa(x_i)}_i(t)=0$. Recall that $\sgn(y)=1$, for $y>0$, and $\sgn(y)=-1$, for $y<0$, whereas $\sgn(0)=0$. Provided that the components of the random walk jumps are independent, broad sufficient conditions for the validity of (\ref{assumptions:large_deviations_tails_right}) are gathered, e.g., in Theorem~5.4.1~(i), (ii) of \cite{Borovkov_Borovkov_08}.

Define the following sets in $\mathbb{R}^d$
\begin{equation}\label{def_O_semiexp}
\mathcal{O}_{\varepsilon}:=\left\{{\bf x}\in\mathbb{R}^d:\sum_{i=1}^d\left|x_i\right|^{\gamma^{\kappa(x_i)}_i}>\nu+\varepsilon\right\},\quad\varepsilon\geq0,\quad\mathcal{O}:=\mathcal{O}_{0},
\end{equation}
\begin{equation}\label{def_Q_semiexp}
\mathcal{Q}_{\varepsilon}:=\left\{{\bf x}\in\mathbb{R}^d:\sum_{i=1}^d\left|x_i\right|^{\gamma^{\kappa(x_i)}_i}<\nu-\varepsilon\right\},\;
\varepsilon\in[0,\nu),\;\mathcal{Q}:=\mathcal{Q}_0,
\end{equation}
\begin{equation}\label{def_P_semiexp}
\mathcal{P}:=\partial\mathcal{O}=\partial\mathcal{Q}=\left\{{\bf x}\in\mathbb{R}^d:\sum_{i=1}^d\left|x_i\right|^{\gamma^{\kappa(x_i)}_i}=\nu\right\},
\end{equation}
where $\partial\mathcal{U}$ stands for the boundary of a set $\mathcal{U}$. Stipulate that ${\bf X}^v(u)/{\bf R}^{-1,\kappa}(t)$ is a vector in $\mathbb{R}^d$ with $i$th coordinate equal to $X^v_i(u)/R^{-1,\kappa(X^v_i(u))}_i(t)$, $u,t\geq0$, $i=1,\ldots,d$.

\begin{Thm}\label{T:main_result_semiexp_lattice}
Let assumptions (\ref{condition1}), (\ref{assumption:tails_right_semiexp})-(\ref{assumptions:large_deviations_tails_right}) be satisfied for supercritical CBRW on $\mathbb{Z}^d$ with Malthusian parameter $\nu$. Then, for each starting point ${\bf z}\in\mathbb{Z}^d$, the following relations are valid.
\begin{equation}\label{T:assertion_1_semiexp}
{\sf P}_{\bf z}\!\left(\omega:\forall\varepsilon>0\;\exists t_1=t_1(\omega,\varepsilon)\;\mbox{s.t.}\;\forall t\geq t_1\;\mbox{and}\;\forall v\in N(t),\;{\bf X}^v(t)/{\bf R}^{-1,\kappa}(t)\notin\mathcal{O}_{\varepsilon}\right)=1,
\end{equation}
\begin{equation}\label{T:assertion_2_semiexp}
{\sf P}_{\bf z}\!\left(\left.\omega:\!\forall\varepsilon\!\in(0,\nu)\,\exists t_2=t_2(\omega,\varepsilon)\;\mbox{s.t.}\;\forall t\geq t_2\;\exists v\in N(t),\;{\bf X}^v(t)/{\bf R}^{-1,\kappa}(t)\notin\mathcal{Q}_{\varepsilon}\right|\mathcal{I}\right)=1,
\end{equation}
where the sets $\mathcal{O}_{\varepsilon}$ and $\mathcal{Q}_{\varepsilon}$ are defined in formulas (\ref{def_O_semiexp}) and (\ref{def_Q_semiexp}).
\end{Thm}

{\bf Remark 1}.
Theorem~\ref{T:main_result_semiexp_lattice} means that, for almost all $\omega$, for any time large enough and any $\varepsilon>0$, there are no particles with properly normalized positions outside the surface $\partial\mathcal{O}_{\varepsilon}$ and, for almost all $\omega\in\mathcal{I}$, there are always such particles outside the surface $\partial\mathcal{Q}_{\varepsilon}$. In other words, for $\omega\in\mathcal{I}$, the most distant particles (``front'' of the population spread) after normalization are located for any time large enough in the interlayer between $\partial\mathcal{O}_{\varepsilon}$ and $\partial\mathcal{Q}_{\varepsilon}$ with $\varepsilon$ small enough. For almost all $\omega\notin\mathcal{I}$, the limit of the normalized particles positions is trivial, that is, equals ${\bf 0}$ (in \cite{B_Doklady_15} one can find necessary and sufficient conditions to guarantee that ${\sf P}(\mathcal{I})=1$). It is natural to call \emph{the limiting shape of the front} the surface $\mathcal{P}=\partial\mathcal{O}=\partial\mathcal{Q}$. Equivalently one can reformulate Theorem~\ref{T:main_result_semiexp_lattice} describing the neighborhoods of $\mathcal{P}$ and $\mathcal{Q}$ in terms of Euclidean distances (instead of taking $\mathcal{O}_{\varepsilon}$ and $\mathcal{Q}_{\varepsilon}$ for $\varepsilon >0$).
\vskip0.2cm
The next result asserts that each point of $\mathcal{P}$ can be considered as a limiting point for the normalized particles positions in CBRW, that is, the surface $\mathcal{P}$ is minimal in a sense.

\begin{Thm}\label{T:one_point_semiexp}
Let conditions of Theorem~\ref{T:main_result_semiexp_lattice} be satisfied. Then, for each ${\bf z}\in\mathbb{Z}^d$ and ${\bf y}\in\mathcal{P}$, one has
$${\sf P}_{\bf z}\left(\left.\omega:\forall t\geq0\;\exists v_{\bf y}=v_{\bf y}(t,\omega)\in N(t)\;\mbox{such that}\;\lim_{t\to\infty}\frac{{\bf X}^{v_{\bf y}}(t)}{{\bf R}^{-1,\kappa}(t)}={\bf y}\right|\mathcal{I}\right)=1.$$
\end{Thm}

{\bf Remark 2}.
Stress that we dot not give a rigorous definition of the front
(but only limiting shape of the front), since it can vary
being dependent on the studied features of the cloud $N(t)$
and basic assumptions. For example, in our framework we could write the definition
of the front as follows. The front of the population propagation is a cloud of particles at time
$t$ such that upon the normalization of particles positions as in
Theorem~\ref{T:main_result_semiexp_lattice} and letting $t\to\infty$,
the almost sure limit points of the particles from the cloud form the surface
$\mathcal{P}$ in (\ref{def_P_semiexp}) called the limiting shape of the front.
However, in other framework such as \cite{B_Arxiv_18} there is a \emph{random limit}
of the properly normalized positions of the most distant (from the origin) particles
in the sense of \emph{weak} convergence. So, the possible mentioned definition of the front
is not suitable for \cite{B_Arxiv_18}. Whenever we talk about the propagation
front for CBRW on $\mathbb{Z}^d$ we mean the generalization to the multidimensional case
of the maximum and the minimum bounding the population on an integer line.
Nevertheless, an attractive possible definition of the front as the set of particles at time $t$
being the most distant on each ray from the origin seems also inconvenient since
continuum of rays will not contain any particle.

\section{Proofs}\label{s:proof_2_semiexp}
To establish Theorems \ref{T:main_result_semiexp_lattice} and \ref{T:one_point_semiexp} we derive a system of non-linear integral equations and estimate its solution, use renewal theory, Laplace transform, coupling, theory of large deviations for the random walk with semi-exponential distributions of jumps, and other probabilistic and analytic technique.
Divide the proof of Theorem~\ref{T:main_result_semiexp_lattice} into 5 Steps.  Within Steps~1,~2, and 3 we consider the case of a single catalyst ${\bf w}_1$ located, without loss of generality, at the origin, that is $W=\{{\bf w}_1\}$ with ${\bf w}_1={\bf 0}$, and the starting point of CBRW being ${\bf 0}$. Within steps 4 and 5 we turn to the general case. In fact, proving
Theorem~\ref{T:main_result_semiexp_lattice} we simultaneously obtain the statement of Theorem~\ref{T:one_point_semiexp}.

\vskip0.2cm
\emph{Step 1.} At this stage we prove statement (\ref{T:assertion_1_semiexp}) in the case of a single catalyst located at ${\bf 0}$ and the starting point ${\bf 0}$ as well.

Let $E(t;\mathcal{U}):={\sf P}_{\bf 0}\left(\exists v\in N(t):{\bf X}^v(t)\in\mathcal{U}\right)$ for a set $\mathcal{U}\subset\mathbb{R}^d$. The following result is a multidimensional counterpart of Lemma~1 in \cite{B_Arxiv_18}
providing an integral equation for the probability $E(t;\mathcal{U})$.

\begin{Lm}\label{L:equation_multi}
Let condition (\ref{condition1}) be valid. Then the probability $E(t;\mathcal{U})$, $t\geq0$, $\mathcal{U}\subset\mathbb{R}^d$, ${\bf 0}\notin\mathcal{U}$, satisfies the non-linear integral equation of convolution type
$$
E(t;\mathcal{U})=\alpha_1\int\nolimits_0^t{ \left(1-f_1\left(1-E(t-s;\mathcal{U})\right)\right)\,dG_1(s)}
$$
\begin{equation}\label{E(t;S)_equation}
+(1-\alpha_1)\int\nolimits_0^t{ E(t-s;\mathcal{U})\,dG_{1,1}(s)}+I\left(t;\mathcal{U}\right),
\end{equation}
where we set
$$I(t;\mathcal{U}):=\sum_{{\bf y}\neq{\bf 0}}{(1-\alpha_1)\frac{q({\bf 0},{\bf y})}{q}\int\nolimits_0^t{ {\sf P}_{\bf y}\left({\bf S}(t-s)\in\mathcal{U},\tau_{{\bf y},{\bf 0}}>t-s\right)\,dG_1(s)}}.$$
\end{Lm}
{\sc Proof.} Similar to the proof of Lemma~1 in \cite{B_Arxiv_18}, consider all the possible evolutions of the parent particle in CBRW. Namely, after time, distributed exponentially with parameter $\beta_1$, it may either produce $k\in\mathbb{Z}_+$ offsprings with probability $\alpha_1{\sf P}(\xi_1=k)$, or jump to the point ${\bf y}\neq{\bf 0}$ with probability $(1-\alpha_1)q({\bf 0},{\bf y})/q$ and first return to the origin in time $\tau_{{\bf y},{\bf 0}}$. If the parent particle does not return to the origin until time $t$, it performs an ordinary random walk ${\bf S}$ starting from ${\bf y}$. At last, it might occur that the parent particle has not undergone changes by time $t$. Summarizing all the above and taking into account (\ref{G-i,G-ij=}), we can write the following formula, for any $\mathcal{U}\subset\mathbb{R}^d$, ${\bf 0}\notin\mathcal{U}$,
\begin{eqnarray*}
1-E(t;\mathcal{U})&=&\alpha_1\sum_{k=0}^{\infty}{\sf P}(\xi_1=k)\int\nolimits_0^t{ \left(1-E(t-s;\mathcal{U})\right)^k\,dG_1(s)}+(1-G_1(t))\\
&+&\sum_{{\bf y}\neq{\bf 0}}{(1-\alpha_1)\frac{q({\bf 0},{\bf y})}{q}\int\nolimits_0^t{ \left(1-E(t-s;\mathcal{U})\right)\,d\left(G_1\ast F_{{\bf y},{\bf 0}}(s)\right)}}\\
&+&\sum_{{\bf y}\neq{\bf 0}}{(1-\alpha_1)\frac{q({\bf 0},{\bf y})}{q}\int\nolimits_0^t{ {\sf P}_{\bf y}\left({\bf S}(t-s)\notin\mathcal{U},\,\tau_{{\bf y},{\bf 0}}>t-s\right)\,dG_1(s)}}.
\end{eqnarray*}
Rewriting the latter equation with respect to unknown function $E(t;\mathcal{U})$ and taking into account the obvious identity
$${\sf P}_{\bf y}\left({\bf S}(s)\in\mathcal{U},\tau_{{\bf y},{\bf 0}}>s\right)=1-F_{{\bf y},{\bf 0}}(s)-{\sf P}_{\bf y}\left({\bf S}(s)\notin\mathcal{U},\tau_{{\bf y},{\bf 0}}>s\right),\quad s\geq0,$$
we get the desired assertion.
Lemma~\ref{L:equation_multi} is proved.

The next lemma provides a convenient form for the function $I$ expressed in terms of the probability ${\sf P}_{\bf 0}\left({\bf S}(t)\in\mathcal{U}\right)$ when ${\bf 0}\notin\mathcal{U}$. Its proof follows the proof of Lemma~2 in \cite{B_Arxiv_18} and is omitted here.

\begin{Lm}\label{L:J-1(t;a)=_semiexp}
Let condition (\ref{condition1}) be satisfied. Then, for any $t\geq0$ and $\mathcal{U}\subset\mathbb{R}^d$, ${\bf 0}\notin\mathcal{U}$, the following identity holds true
\begin{eqnarray}\label{J-1(t;a)=_semiexp}
\frac{q I(t;\mathcal{U})}{(1-\alpha_1)\beta_1}&=&{\sf P}_{\bf 0}\left({\bf S}(t)\in\mathcal{U}\right)-\int\nolimits_0^t{ {\sf P}_{\bf 0}\left({\bf S}(t-s)\in\mathcal{U}\right)\,d F_{{\bf 0},{\bf 0}}(s)}\\
&-&\frac{\beta_1-q}{\beta_1}\int\nolimits_0^t{ {\sf P}_{\bf 0}\left({\bf S}(t-s)\in\mathcal{U}\right)\,d\left(G_1(s)-G_1\ast F_{{\bf 0},{\bf 0}}(s)\right)}.\nonumber
\end{eqnarray}
\end{Lm}

The definition of supercritical regime of CBRW implies that $\alpha_1m_1+(1-\alpha_1)F_{{\bf 0},{\bf 0}}(\infty)>1$
and $\alpha_1m_1G_1^{\ast}(\nu)+(1-\alpha_1)\,G_1^{\ast}(\nu)\overline{F}^{\,\ast}_{{\bf 0},{\bf 0}}(\nu)=1$.
In terms of the function $G(t):=\alpha_1m_1G_1(t)+(1-\alpha_1)\,G_1\ast\overline{F}_{{\bf 0},{\bf 0}}(t)$, $t\geq0$, it means that $G^{\ast}(\nu)=1$.

\begin{Lm}\label{L:E(t;)_estimate_semiexp}
Let conditions (\ref{condition1}), (\ref{assumption:tails_right_semiexp}), and (\ref{assumptions:large_deviations_tails_right}) be satisfied. Fix ${\bf x}=(x_1,\ldots,x_d)$ from the set $\partial\mathcal{O}_{\varepsilon}\cap\left\{{\bf x}\in\mathbb{R}^d:x_i\geq0,i=1,\ldots,d\right\}=\left\{{\bf x}\in\mathbb{R}^d_+:\sum_{i=1}^d x_i^{\gamma^+_i}=\nu+\varepsilon\right\}=:\partial\mathcal{O}^+_{\varepsilon}$. Then there exists $\varepsilon_0=\varepsilon_0(\nu,\varepsilon)>0$ such that
\begin{equation}\label{E(t;U)_estimate_semiexp}
E(t;\Delta({\bf x};t))\leq C e^{-\varepsilon_0t},\quad t\geq0,
\end{equation}
where $\Delta({\bf x};t):=\left[x_1R^{-1,+}_1(t),+\infty\right)\times\ldots\times\left[x_dR^{-1,+}_d(t),+\infty\right)\subset\mathbb{R}^d$ and $C$ is a positive constant.
\end{Lm}
{\sc Proof.} For any $\mathcal{U}\subset\mathbb{R}^d$, ${\bf 0}\notin\mathcal{U}$, by mean value theorem on $f_1$, equation (\ref{E(t;S)_equation}) entails the inequality
$$E\left(t;\mathcal{U}\right)\leq\int\nolimits_0^t{ E\left(t-s;\mathcal{U}\right)\,d G(s)}+I\left(t;\mathcal{U}\right).$$
Iterating this inequality $k$ times we get
$$E\left(t;\mathcal{U}\right)\leq\int\nolimits_0^t{ E\left(t-s;\mathcal{U}\right)\,d G^{\ast (k+1)}(s)}+\int\nolimits_0^t{ I\left(t-s;\mathcal{U}\right)\,d\sum_{j=0}^k{G^{\ast j}(s)}}.$$
For any fixed $t$, one has $G^{\ast k}(t)\to 0$, as $k\to\infty$, for example, due to Lemma~22 in \cite{Vatutin_book_09}. Hence, the term $\int\nolimits_0^t{ E\left(t-s;\mathcal{U}\right)\,dG^{\ast(k+1)}(s)}$ is negligibly small for large $k$. The latter inequality can be rewritten as follows
\begin{equation}\label{E(t;U)_inequality_semiexp}
E\left(t;\mathcal{U}\right)\leq\int\nolimits_0^t{ I\left(t-s;\mathcal{U}\right)\,d\sum_{j=0}^{\infty}{G^{\ast j}(s)}}.
\end{equation}

It follows from (\ref{J-1(t;a)=_semiexp}) that
$$I\left(t;\mathcal{U}\right)\leq\frac{(1-\alpha_1)\beta_1}{q}{\sf P}_{\bf 0}\left({\bf S}(t)\in\mathcal{U}\right)+\frac{(1-\alpha_1)\left|\beta_1-q\right|}{q}\int\nolimits_0^t{ {\sf P}_{\bf 0}\left({\bf S}(t-s)\in\mathcal{U}\right)\,dG_1(s)}.$$
Combining this inequality with (\ref{E(t;U)_inequality_semiexp}) we get
$$E\left(t;\mathcal{U}\right)\leq\int\nolimits_{0}^{t}{ {\sf P}_{\bf 0}({\bf S}(t-s)\in\mathcal{U})\,d\!\left(\frac{(1-\alpha_1)\beta_1}{q}\sum_{k=0}^{\infty}G^{\ast k}(s)\!+\!\frac{(1-\alpha_1)\left|\beta_1-q\right|}{q}G_1\ast\sum_{k=0}^{\infty}G^{\ast k}(s)\right)}\!.$$
Consider $\mathcal{U}=\Delta({\bf x};t)$, where ${\bf x}\in\partial\mathcal{O}^+_{\varepsilon}$. Then by virtue of assumptions (\ref{assumption:tails_right_semiexp}), (\ref{assumptions:large_deviations_tails_right}), and Theorem~25 in \cite{Vatutin_book_09}, p.~30, for any $\delta_1$ and $\delta_2$ from $(0,1)$, there exists $T=T(\delta_1,\delta_2)$ such that, for any $t\geq T$, one has
$$E\left(t;\Delta\left({\bf x};t\right)\right)\leq C_1\int\nolimits_{0}^{t}{ h(t-s)\,d\sum_{k=0}^{\infty}G^{\ast k}(s)}\;\prod_{i=1}^d{\sf P}\left(Y_i^1\geq x_iR^{-1,+}_i(t)\right)^{1-\delta_1}$$
$$\leq C_2\,e^{\nu t}\frac{\int\nolimits_0^{ \infty}h(s)e^{-\nu s}\,ds}{\int\nolimits_0^{\infty}{se^{-\nu s}\,d G(s)}}\prod_{i=1}^d{R^+_i\left(x_iR^{-1,+}_i(t)\right)^{1-\delta_1}}\leq C_3\,e^{\nu t}\prod_{i=1}^d{\exp\left\{-(1-\delta_1)x^{\gamma^+_i}_i t(1+o(1))\right\}}$$
$$\leq C_3\,\exp\left\{\nu t-(1-\delta_1)t\sum_{i=1}^d{x^{\gamma^+_i}_i}+\delta_2t\right\}\leq C_3\exp\left\{-\left((1-\delta_1)(\nu+\varepsilon)-\nu-\delta_2\right)t\right\},$$
for some positive constants $C_1$, $C_2$, and $C_3$. One can take $\delta_1,\delta_2\in(0,1)$ in such a way that $(1-\delta_1)(\nu+\varepsilon)-\nu-\delta_2=\varepsilon_0>0$. Lemma~\ref{L:E(t;)_estimate_semiexp} is proved.

\begin{Lm}\label{L:T:assertion_1_semiexp}
Let conditions (\ref{condition1}), (\ref{assumption:tails_right_semiexp}) and (\ref{assumptions:large_deviations_tails_right}) be valid. Then the following relation holds true
\begin{equation}\label{Step1:T:assertion_1_semiexp}
{\sf P}_{\bf 0}\!\left(\omega:\forall\varepsilon>0\;\exists t_3=t_3(\omega,\varepsilon)\;\mbox{s.t.}\;\forall t\geq t_3\;\mbox{and}\;\forall v\in N(t),\;{\bf X}^v(t)/{\bf R}^{-1,\kappa}(t)\notin\mathcal{O}_{\varepsilon}\cap\mathbb{R}^d_+\right)=1.
\end{equation}
\end{Lm}
{\sc Proof.} Fix $j\in\mathbb{N}$ and ${\bf x}\in\partial\mathcal{O}^+_{\varepsilon+1/j}$.
Set $A_t:=\{\omega:\forall v\in N(t)\;\mbox{one
has}\;{\bf X}^v(t)\notin\Delta({\bf x};t)\}$, $t\geq 0$.
As usual, $\overline{A}$ stands for the complement of a set $A$ and $\{A_n\;\mbox{infinitely often}\;\}=\{A_n\;\mbox{i.o.}\}=\cap_{k=1}^{\infty}\cup_{n=k}^{\infty}A_n$, for a sequence of sets $A_n$. By virtue of Borel-Cantelli lemma the estimate (\ref{E(t;U)_estimate_semiexp}) entails ${\sf P}_{\bf
0}\left(\overline{A}_{n/2^m}\;\mbox{i.o.}\right)=0$, for any fixed
$m\in\mathbb{N}$. Consequently, ${\sf P}_{\bf
0}\left(\cap_{m=1}^{\infty}\cup_{k=1}^{\infty}\cap_{n=k}^{\infty}A_{n/2^m}\right)=1$.
It means that, for almost all $\omega\in\Omega$ and for any
$m\in\mathbb{N}$, there exists a positive integer $k_1=k_1(m,\omega)$ such
that, for any $n\geq k_1$ and every $v\in N(n/2^m)$, one has
${\bf X}^v(n/2^m)\notin\Delta({\bf x};n/2^m)$. Since the
set of binary rational numbers is dense in $\mathbb{R}$ and the
sojourn time of a particle $v\in N(t)$ in a set $\Delta({\bf x};t)$ (conditioned on the event that the particle has hit the set) contains non-zero interval with probability $1$, we
conclude that
\begin{equation}\label{P(O_gamma,epsilon)=1_semiexp}
{\sf P}_{\bf 0}\left(\omega:\exists t_4(\omega)\;\mbox{such that}\;\forall t\geq t_4(\omega)\,\mbox{and}\,\forall v\in N(t),\;{\bf X}^v(t)\notin\Delta({\bf x};t)\right)=1.
\end{equation}

Unfix ${\bf x}\in\partial\mathcal{O}^+_{\varepsilon+1/j}$. If the set $\partial\mathcal{O}^+_{\varepsilon+1/j}$ is finite (it occurs when $d=1$), put $\Upsilon_j=\partial\mathcal{O}^+_{\varepsilon+1/j}$. Otherwise, let $\Upsilon_j$ be an everywhere dense set in $\partial\mathcal{O}^+_{\varepsilon+1/j}$ joined with points ${\bf x}\in\mathbb{R}^d_+$ with $x_i=0$, $i=1,\ldots,d$, $i\neq l$, for each $l=1,\ldots,d$. For instance, let $\Upsilon_j$ be the set of vectors ${\bf x}$ from $\partial\mathcal{O}^+_{\varepsilon+1/j}$ with rational coordinates $x_i$, $i=1,\ldots,d$, $i\neq l$, for each $l=1,\ldots,d$.

Unfix $j\in\mathbb{N}$. Consider the domain $\mathcal{O}_{\varepsilon}\cap\mathbb{R}^d_+=\cup_{j=1}^{\infty}\cup_{{\bf x}\in\Upsilon_j}[x_1,+\infty)\times\ldots\times[x_d,+\infty)$. Take into account that the relation ${\bf X}^v(t)\notin\Delta({\bf x};t)$ is equivalent to ${\bf X}^v(t)/{\bf R}^{-1,\kappa}(t)\notin[x_1,+\infty)\times\ldots\times[x_d,+\infty)$. Then formula \eqref{P(O_gamma,epsilon)=1_semiexp} entails the equality
$${\sf P}_{\bf 0}\left(\omega:\exists t_5(\omega)\;\mbox{such that}\;\forall t\geq t_5(\omega)\,\mbox{and}\,\forall v\in N(t),\;{\bf X}^v(t)/{\bf R}^{-1,\kappa}(t)\notin\mathcal{O}_{\varepsilon}\cap\mathbb{R}^d_+\right)=1,$$
valid for each $\varepsilon>0$. Hence the latter relation implies the assertion (\ref{Step1:T:assertion_1_semiexp}). Lemma~\ref{L:T:assertion_1_semiexp} is proved.

When $d=1$, Lemma~\ref{L:T:assertion_1_semiexp} states that $\limsup\nolimits_{t\to\infty}M_t/(t^{1/\gamma^+_1}L^{(3,+)}_1(t))\leq\nu^{1/\gamma^+_1}$ a.s. Thus, we obtain the upper estimate for the maximum $M_t$ in the case of CBRW with a single catalyst at $0$ and the starting point $0$.

Let $d\in\mathbb{N}$. In Lemma~\ref{L:T:assertion_1_semiexp} we consider the particles propagation in the positive orthant $\mathbb{R}^d_+$. Now trace the spread of particles with growing time in other directions. Without loss of generality, we deal with $\mathbb{R}_{-}\times\mathbb{R}^{d-1}_+$. Reflect the lattice $\mathbb{Z}^d$ with particles in CBRW on it at each time $t$ with respect to plane $x_1=0$. We get a new CBRW on $\mathbb{Z}^d$ and may apply
to it Lemma~\ref{L:T:assertion_1_semiexp}. Consequently,
\begin{equation}\label{reformulation_Lemma_semiexp}
{\sf P}_{\bf 0}\!\left(\omega:\forall\varepsilon>0\;\exists t_6=t_6(\omega,\varepsilon)\;\mbox{s.t.}\;\forall t\geq t_6\;\mbox{and}\;\forall v\in N(t),{\bf X}^v(t)/{\bf R}^{-1,\kappa}(t)\notin\mathcal{O}_{\varepsilon}\!\cap\!\left(\mathbb{R}_{-}\!\times\!\mathbb{R}^{d-1}_+\right)\right)\!=\!1.
\end{equation}
In the same manner reformulation of Lemma~\ref{L:T:assertion_1_semiexp} for other orthants in $\mathbb{R}^d$ combined with (\ref{Step1:T:assertion_1_semiexp}) and (\ref{reformulation_Lemma_semiexp}) leads to the first assertion of Theorem~\ref{T:main_result_semiexp_lattice} in the case of CBRW with a single catalyst at ${\bf 0}$ and the starting point ${\bf 0}$.

\vskip0.2cm
\emph{Step 2.}
We prove statement (\ref{T:assertion_2_semiexp}) whenever there is a single catalyst located at ${\bf 0}$ and the starting point is ${\bf 0}$ as well. We temporarily assume that ${\sf E}\xi^2_1<\infty$ and follow the ideas of \cite{Carmona_Hu_14}, Sect.5.2.

\begin{Lm}\label{L:lower_estimate_semiexp}
Let conditions (\ref{condition1}), (\ref{assumption:tails_right_semiexp}), and (\ref{assumptions:large_deviations_tails_right}) be satisfied. Choose function $r=r(t)$ in such a way that $r(t)\leq t$, $r(t)\nearrow+\infty$, $t-r(t)\nearrow+\infty$ and $t-r(t)=o(t)$, as $t\to\infty$ (for example, we can put $r(t)=t-\ln{t}$). Fix both $\varepsilon\in(0,\nu)$ and ${\bf x}\in\partial\mathcal{Q}^+_{\varepsilon}:=\partial\mathcal{Q}_{\varepsilon}\cap\mathbb{R}^d_+$. Then, for some positive constant $C_4$, one has
\begin{equation}\label{E(t;U)_lower_estimate_semiexp}
{\sf P}_{\bf 0}\left({\bf X}^v(t)\notin\Delta({\bf x};t)\mbox{\;for any\;}v\in N(t),\mu(r;{\bf 0})\geq C_4 e^{\nu r}\right)\leq\exp\left\{-e^{\varepsilon t+o(t)}\right\},\quad t\to\infty.
\end{equation}
\end{Lm}
{\sc Proof.} In view of Theorem~4 in \cite{B_Doklady_15}, on the set $\mathcal{I}$ at time $r$, $0<r<t$, there are at least $[C_4e^{\nu r}]$ particles at ${\bf 0}$ for some positive constant $C_4$ (as usual, $[r]$ stands for the integer part of a number $r\in\mathbb{R}_+$). If these particles move according to the random walk ${\bf S}$ such that ${\bf S}(u)\neq{\bf 0}$ for each $u\in[\tau_{\bf 0},t-r]$, then remote particles in CBRW at time $t$ are not less far than $[C_4e^{\nu r}]$ i.i.d. copies of ${\bf S}(t-r)$ with ${\bf S}(u)\neq{\bf 0}$, for each $u\in[\tau_{\bf 0},t-r]$.

For a set $\mathcal{U}\subset\mathbb{R}^d$, ${\bf 0}\notin\mathcal{U}$, and $t\geq0$, the following identity is valid
$${\sf P}_{\bf 0}\left({\bf S}(t)\in\mathcal{U},\tau_{{\bf 0},{\bf 0}}>t\right)={\sf P}_{\bf 0}\left({\bf S}(t)\in\mathcal{U}\right)-\int\nolimits_0^t{ {\sf P}_{\bf 0}\left({\bf S}(t-s)\in\mathcal{U}\right)\,d F_{{\bf 0},{\bf 0}}(s)}.$$
Then assumptions~(\ref{assumption:tails_right_semiexp}) and (\ref{assumptions:large_deviations_tails_right}) imply that
$${\sf P}_{\bf 0}\left({\bf S}(t-r)\in\Delta({\bf x};t),\tau_{{\bf 0},{\bf 0}}>t-r\right)$$
$$={\sf P}_{\bf 0}\left({\bf S}(t-r)\in\Delta({\bf x};t)\right)-\int\nolimits_0^{t-r}{ {\sf P}_{\bf 0}\left({\bf S}(t-r-s)\in\Delta({\bf x};t)\right)\,d F_{{\bf 0},{\bf 0}}(s)}$$
$$=\left(h(t-r)-\int\nolimits_0^{t-r}{ h(t-r-s)\,dF_{{\bf 0},{\bf 0}}(s)}+o\left(h(t-r)\right)\right)\prod_{i=1}^d{\sf P}\left(Y_i^1\geq x_iR^{-1,+}_i(t)\right)^{1+\varepsilon^{+}_i(t)}$$
$$=\left(h(t-r)\left(1-F_{{\bf 0},{\bf 0}}(t-r)\right)+\int_0^{t-r}{ \left(h(t-r)-h(t-r-s)\right)\,dF_{{\bf 0},{\bf 0}}(s)}+o\left(h(t-r)\right)\right)$$
$$\times\prod_{i=1}^d{R^+_i\left(x_iR^{-1,+}_i(t)\right)^{1+\varepsilon^+_i(t)}}=\exp\left\{-t\sum_{i=1}^d{x_i^{\gamma_i^{+}}}+o(t)\right\}=\exp\left\{-(\nu-\varepsilon)t+o(t)\right\},\quad t\to\infty.$$
The latter relation leads to the estimate
$${\sf P}_{\bf 0}\left({\bf X}^v(t)\notin\Delta({\bf x};t)\mbox{\;for any\;}v\in N(t),\mu(r;{\bf 0})\geq C_4 e^{\nu r}\right)$$
$$\leq\left(1-{\sf P}_{\bf 0}\left({\bf S}(t-r)\in\Delta({\bf x};t),\tau_{{\bf 0},{\bf 0}}>t-r\right)\right)^{\left[C_4 e^{\nu r}\right]}$$
$$\leq\exp\left\{-\left[C e^{\nu r}\right]e^{-(\nu-\varepsilon)t+o(t)}\right\}=\exp\left\{-e^{\nu r-(\nu-\varepsilon)t+o(t)}\right\},\quad t\to\infty.$$
The assertion of Lemma~\ref{L:lower_estimate_semiexp} now follows from our choice of $r=r(t)$. Lemma~\ref{L:lower_estimate_semiexp} is proved.

\begin{Lm}\label{L:T:assertion_2_semiexp}
Let conditions (\ref{condition1}), (\ref{assumption:tails_right_semiexp}), and (\ref{assumptions:large_deviations_tails_right}) be valid. Then the following relation holds true
\begin{equation}\label{Step2:T:assertion_2_semiexp}
{\sf P}_{\bf 0}\!\left(\left.\omega:\forall\varepsilon\in(0,\nu)\;\exists t_7=t_7(\omega,\varepsilon)\;\mbox{s.t.}\;\forall t\geq t_7\;\exists v\in N(t),\;{\bf X}^v(t)/{\bf R}^{-1,\kappa}(t)\notin\mathcal{Q}_{\varepsilon}\cap\mathbb{R}^d_+\right|\mathcal{I}\right)=1.
\end{equation}
\end{Lm}
{\sc Proof.} Fix ${\bf x}\in\partial\mathcal{Q}^+_{\varepsilon}$.
Denote by $B_t$ the event $\{\omega:\exists v\in N(t)\;\mbox{such that}\;{\bf X}^v(t)\in\Delta({\bf x};t)\}$. By virtue of Borel-Cantelli lemma and Theorem~4 in \cite{B_Doklady_15} the estimate (\ref{E(t;U)_lower_estimate_semiexp}) entails ${\sf P}_{\bf
0}\left(\left.\overline{B}_{n/2^m}\;\mbox{i.o.}\right|\mathcal{I}\right)=0$, for any fixed
$m\in\mathbb{N}$. Consequently, ${\sf P}_{\bf
0}\left(\left.\cap_{m=1}^{\infty}\cup_{k=1}^{\infty}\cap_{n=k}^{\infty}B_{n/2^m}\right|\mathcal{I}\right)=1$.
It means that for almost all $\omega\in\Omega$ and for any
$m\in\mathbb{N}$ there exists positive integer $k_2=k_2(m,\omega)$ such
that, for any $n\geq k_2$, one can find $v\in N(n/2^m)$ such that
${\bf X}^v(n/2^m)\in\Delta({\bf x};n/2^m)$. Since the
set of binary rational numbers is dense in $\mathbb{R}$ and the
sojourn time of a particle $v\in N(t)$ in a set $\Delta({\bf x};t)$ (conditioned on the event that the particle has hit the set) contains non-zero interval with probability $1$, we
conclude that
\begin{equation}\label{P(O_gamma,epsilon)=1_semiexp_lower}
{\sf P}_{\bf 0}\left(\left.\omega:\exists t_8(\omega)\;\mbox{such that}\;\forall t\geq t_8(\omega)\,\mbox{one has}\,\exists v\in N(t),\;{\bf X}^v(t)\in\Delta({\bf x};t)\right|\mathcal{I}\right)=1.
\end{equation}

Unfix ${\bf x}\in\partial\mathcal{Q}^+_{\varepsilon}$. If the set $\partial\mathcal{Q}^+_{\varepsilon}$ is finite (it occurs when $d=1$), put $\Upsilon=\partial\mathcal{Q}^+_{\varepsilon}$. Otherwise, let $\Upsilon$ be an everywhere dense set in $\partial\mathcal{Q}^+_{\varepsilon}$ joined with points ${\bf x}\in\mathbb{R}^d_+$ with $x_i=0$, $i=1,\ldots,d$, $i\neq l$, for each $l=1,\ldots,d$. For instance, let $\Upsilon$ be the set of vectors ${\bf x}$ from $\partial\mathcal{Q}^+_{\varepsilon}$ with rational coordinates $x_i$, $i=1,\ldots,d$, $i\neq l$, for each $l=1,\ldots,d$.

Consider the domain $\mathcal{Q}_{\varepsilon}\cap\mathbb{R}^d_+=\mathbb{R}^d_+\setminus\cup_{{\bf x}\in\Upsilon}[x_1,+\infty)\times\ldots\times[x_d,+\infty)$. Take into account that the relation ${\bf X}^v(t)\in\Delta({\bf x};t)$ is equivalent to ${\bf X}^v(t)/{\bf R}^{-1,\kappa}(t)\in[x_1,+\infty)\times\ldots\times[x_d,+\infty)$. Then formula \eqref{P(O_gamma,epsilon)=1_semiexp_lower} entails the equality
$${\sf P}_{\bf 0}\left(\left.\omega:\exists t_9(\omega)\;\mbox{such that}\;\forall t\geq t_9(\omega)\,\exists v\in N(t),\;{\bf X}^v(t)/{\bf R}^{-1,\kappa}(t)\notin\mathcal{Q}_{\varepsilon}\cap\mathbb{R}^d_+\right|\mathcal{I}\right)=1,$$
valid for each $\varepsilon\in(0,\nu)$. Unfix $\varepsilon\in(0,\nu)$. Then the latter relation implies the assertion (\ref{Step2:T:assertion_2_semiexp}). Lemma~\ref{L:T:assertion_2_semiexp} is proved.

When $d=1$, Lemma~\ref{L:T:assertion_2_semiexp} states that $\liminf\nolimits_{t\to\infty}M_t/(t^{1/\gamma^+_1}L^{(3,+)}_1(t))\geq\nu^{1/\gamma^+_1}$ a.s. on the set $\mathcal{I}$. Thus, we obtain the lower estimate for the maximum $M_t$ in the case of CBRW with a single catalyst at $0$ and the starting point $0$ under the additional assumption ${\sf E}\xi^2_1<\infty$.

Let $d\in\mathbb{N}$. In Lemma~\ref{L:T:assertion_2_semiexp} we consider the particles propagation in the positive orthant $\mathbb{R}^d_+$. Similarly to discussion at the end of \emph{Step~1}, we may reformulate Lemma~\ref{L:T:assertion_2_semiexp} for other orthants in $\mathbb{R}^d$. Reformulation of Lemma~\ref{L:T:assertion_2_semiexp} for other orthants in $\mathbb{R}^d$ combined with (\ref{Step2:T:assertion_2_semiexp}) leads to the second assertion of Theorem~\ref{T:main_result_semiexp_lattice} in the case of CBRW with a single catalyst at ${\bf 0}$ and the starting point ${\bf 0}$ whenever ${\sf E}\xi^2_1<\infty$.

Combination of the proved in \emph{Step 1} assertion (\ref{T:assertion_1_semiexp}) and relation (\ref{P(O_gamma,epsilon)=1_semiexp_lower}), valid for each ${\bf x}\in\partial\mathcal{Q}_{\varepsilon}$, implies the statement of Theorem~\ref{T:one_point_semiexp} for the case of a single catalyst at ${\bf 0}$ and the starting point ${\bf 0}$ whenever ${\sf E}\xi^2_1<\infty$.

\vskip0.2cm \emph{Step 3.}
Assume that $W=\{{\bf w}_1\}$ with ${\bf w}_1={\bf 0}$ and the starting point of CBRW is
${\bf 0}$ whereas now ${\sf E}\xi^2_1=\infty$. To verify assertion (\ref{P(O_gamma,epsilon)=1_semiexp_lower}) (and, as a consequence, (\ref{T:assertion_2_semiexp})) under such assumptions one can follow the
proof scheme proposed in \cite{Carmona_Hu_14}, Sec.~5.3, based on
a coupling. In contrast to \cite{Carmona_Hu_14} we employ Theorem~3 of \cite{B_Doklady_15}
devoted to the strong convergence of the total and local particles
numbers in supercritical CBRW instead of using properties of a
fundamental martingale as in \cite{Carmona_Hu_14}. Moreover, here
we exploit function $g(u)=\alpha_1
f_1\left(q_{esc}+(1-q_{esc})u\right)+(1-\alpha_1)q_{esc}-u$,
${u\in[0,1]}$, where $q_{esc}={\sf P}_{\bf 0}\left(\overline{\tau}_{{\bf 0},{\bf 0}}=\infty\right)=1-\overline{F}_{{\bf 0},{\bf 0}}(\infty)$ is the
escape probability of the random walk ${\bf S}$. Other details of
the Step~3 proof can be omitted.

\vskip0.2cm \emph{Step 4.} Now we deal with $N>1$ and ${\bf x}\in W$, say ${\bf x}={\bf w}_i$. Let us discuss here the main differences between the case of single and multiple catalysts and sketch the subsequent proof omitting cumbersome details. In the multiple setting the counterpart of the probability $E(t;\mathcal{U})$ is the vector ${\bf E}(t;\mathcal{U}):=\left(E_{{\bf w}_1}(t;\mathcal{U}),\ldots,E_{{\bf w}_N}(t;\mathcal{U})\right)$ with $E_{{\bf w}_i}(t;\mathcal{U}):={\sf P}_{{\bf w}_i}\left(\exists v\in N(t):{\bf X}^v(t)\in\mathcal{U}\right)$, $i=1,\ldots,N$, $t\geq0$, for $\mathcal{U}\subset\mathbb{R}^d$. Similarly to Lemma~\ref{L:equation_multi}, it satisfies the following system of non-linear integral equations of convolution type
\begin{eqnarray*}
E_{{\bf w}_i}(t;\mathcal{U})&=&\alpha_i\int\nolimits_0^t{ \left(1-f_i\left(1-E_{{\bf w}_i}(t-s;\mathcal{U})\right)\right)\,dG_i(s)}\\
&+&(1-\alpha_i)\sum_{j=1}^N\int\nolimits_0^t{ E_{{\bf w}_j}(t-s;\mathcal{U})\,dG_{i,j}(s)}+I_{{\bf w}_i}\left(t;\mathcal{U}\right),
\end{eqnarray*}
where
$$I_{{\bf w}_i}(t;\mathcal{U}):=\sum_{{\bf y}\notin W}{(1-\alpha_i)\frac{q({\bf w}_i,{\bf y})}{q}\int\nolimits_0^t{ {\sf P}_{\bf y}\left({\bf S}(t-s)\in\mathcal{U},{_{W_k}\tau_{{\bf y},{\bf w}_k}}>t-s,k=1,\ldots,N\right)\,dG_i(s)}},$$
$t\geq0$, $i=1,\ldots,N$, $\mathcal{U}\subset\mathbb{R}^d$, $W\cap\mathcal{U}=\varnothing$.

A multiple setting counterpart of Lemma~\ref{L:J-1(t;a)=_semiexp} states that
\begin{eqnarray*}
& &\frac{q I_{{\bf w}_i}(t;\mathcal{U})}{(1-\alpha_i)\beta_i}={\sf P}_{{\bf w}_i}\left({\bf S}(t)\in\mathcal{U}\right)-\sum_{k=1}^N\int\nolimits_0^t{ {\sf P}_{{\bf w}_k}\left({\bf S}(t-s)\in\mathcal{U}\right)\,d\,{_{W_k}F_{{\bf w}_i,{\bf w}_k}(s)}}\\
&-&\frac{\beta_i-q}{\beta_i}\int\nolimits_0^t{ {\sf P}_{{\bf w}_i}\left({\bf S}(t-s)\in\mathcal{U}\right)\,d G_i(s)}+\sum_{k=1}^N\frac{\beta_i-q}{\beta_i}\int\nolimits_0^t{ {\sf P}_{{\bf w}_k}\left({\bf S}(t-s)\in\mathcal{U}\right)\,d G_i\ast{_{W_k}F_{{\bf w}_i,{\bf w}_k}(s)}}.
\end{eqnarray*}
It follows that
$$I_{{\bf w}_i}(t;\mathcal{U})\leq\frac{(1-\alpha_i)\beta_i}{q}{\sf P}_{{\bf w}_i}\left({\bf S}(t)\in\mathcal{U}\right)+\frac{(1-\alpha_i)\left|\beta_i-q\right|}{q}\int\nolimits_0^t{ {\sf P}_{{\bf w}_i}\left({\bf S}(t-s)\in\mathcal{U}\right)\,d G_i(s)}=:K_i(t;\mathcal{U}).$$

While proving Lemma~\ref{L:E(t;)_estimate_semiexp} we obtained inequality (\ref{E(t;U)_inequality_semiexp}). Likewise we derive the following vector inequality, valid coordinate-wise,
$${\bf E}(t;\mathcal{U})\leq\sum_{k=0}^{\infty}{\bf G}^{\ast k}\ast{\bf K}(t;\mathcal{U}).$$
Here a multiple setting counterpart of function $G$ is a matrix ${\bf G}(t)=\left(G^{\,(N)}_{i,j}(t)\right)_{i,j=1}^N$ with entries $G^{\,(N)}_{i,j}(t):=\delta_{i,j}\alpha_i m_i G_i(t)+(1-\alpha_i)G_i\ast{_{W_j}\overline{F}_{{\bf w}_i,{\bf w}_j}(t)}$, $t\geq0$, whereas ${\bf K}(t;\mathcal{U})$, $t\geq0$, is the vector-column with $i$th coordinate $K_i(t;\mathcal{U})$. The element $d_{i,j}(\lambda)$ of matrix $D(\lambda)$, $\lambda\geq0$, is just the Laplace transform of $G^{\,(N)}_{i,j}$.
Recall that the operation ``$\ast$'' of convolution of matrices is defined exactly as matrix multiplication except that we convolve elements rather than multiply them.

As for validating Lemma~\ref{L:E(t;)_estimate_semiexp} for $N=1$, in case $N>1$ we inspect the asymptotic behavior of $\sum_{k=0}^{\infty}{\bf G}^{\ast k}\ast{\bf K}(t;\mathcal{U})$ when $\mathcal{U}=\Delta({\bf x};t)$ with ${\bf x}\in\partial\mathcal{O}^+_{\varepsilon}$ and $t\to\infty$. Employing Corollary~3.1, item (i), in \cite{Crump_70}, we deduce the same estimate as (\ref{E(t;U)_estimate_semiexp}) after replacing $E(t;\Delta({\bf x};t))$ by $E_{{\bf w}_i}(t;\Delta({\bf x};t))$. The rest of the proofs of Theorems~\ref{T:main_result_semiexp_lattice} and \ref{T:one_point_semiexp} in case of CBRW with general catalysts set $W$ and the starting point from $W$ is implemented similar to the arguments of Steps~1 -- 3.

\vskip0.2cm {\it Step 5.} Turning to a supercritical CBRW on $\mathbb{Z}^d$ with a finite catalysts set $W$ and the starting
point ${\bf z}\notin W$, we supplement the catalysts set $W$ with ${\bf w}_{N+1}={\bf x}$ and put $\alpha_{N+1}=0$, $m_{N+1}=0$,
$G_{N+1}(t)=1-e^{-qt}$, $t\geq0$. According to Lemma~3 in \cite{B_TPA_15} a new CBRW with catalysts set $\{{\bf w}_1,\ldots,{\bf w}_{N+1}\}$ is supercritical whenever the underlying CBRW is supercritical and the Malthusian parameters in these CBRW coincide. Then one can apply the proved part of Theorems~\ref{T:main_result_semiexp_lattice} and \ref{T:one_point_semiexp} to the new CBRW and obtain the desired assertions of those theorems for CBRW with an arbitrary starting point.

\vskip0.2cm The proof of Theorems~\ref{T:main_result_semiexp_lattice} and \ref{T:one_point_semiexp} is complete.

\vskip0.2cm
{\bf Remark 3}.
Within the proofs we track the evolution of the particles ``at the front''. As for CBRW with regularly varying tails in \cite{B_Arxiv_18}, a particle ``at the front'' at time $t$ was born at time $t-o(t)$ and then reached the front within time $o(t)$.
This significantly differs from the case of ``light'' tails.
It follows from the proof of Theorem~1.1 in \cite{Carmona_Hu_14}
and Theorem~1 in \cite{B_SPA_18} treating CBRW with light tails that
a particle ``at the front'' at time $t$ was born at time $\theta t$ for
specified $\theta\in(0,1)$ and then walked only until time $t$.

\section{Examples}\label{s:examples_semiexp}

According to formula (\ref{S(t)=representation}) the random walk ${\bf S}$ is a jump process with increments ${\bf Y}^j$. In this section we assume that the coordinates of each jump ${\bf Y}^j=\left(Y^j_1,\ldots,Y^j_d\right)$, $j\in\mathbb{N}$, are independent. Without loss of generality, consider ${\bf x}\in\mathbb{R}^d_+$, ${\bf x}\neq{\bf 0}$. Then
\begin{eqnarray}\label{example1}
{\sf P}_{\bf 0}\left({\bf S}(u)\in\Delta({\bf x};t)\right)&=&\sum_{k=1}^{\infty}{\sf P}\left(\Pi(u)=k\right){\sf P}\left(\sum_{j=1}^k{{\bf Y}^j}\in\Delta({\bf x};t)\right)\nonumber\\
&=&\sum_{k=1}^{\infty}{\sf P}\left(\Pi(u)=k\right)\prod_{i=1}^d{\sf P}\left(\sum_{j=1}^kY^j_i\geq x_iR^{-1,+}_i(t)\right).
\end{eqnarray}

Write $Y=Y^+-Y^-$, where $Y^+:=Y\mathbb{I}\{Y\geq0\}$ and $Y^-:=-Y\mathbb{I}\{Y<0\}$.
Consider $Y^j_i=Y^{(j,+)}_i-Y^{(j,-)}_i$, where $Y^{(j,+)}_i$, $Y^{(j,-)}_i\geq0$ are defined in the mentioned way. Set ${\sf P}\left(Y^{(j,\kappa)}_i=0\right)=1-L^{(1,\kappa)}_i$, $\kappa\in\{+,-\}$, $j=1,2,\ldots$. Then $L^{(1,+)}_i+L^{(1,-)}_i=1$,  $i=1,\ldots,d$. Let also $Y^{(j,\kappa)}_i$ conditioned to be strictly positive have a discrete Weibull distribution, \cite{FKZ_11}, p.~10, with parameters $\gamma^{\kappa}_i$ and $\left(L^{(2,\kappa)}_i\right)^{-1/\gamma^{\kappa}_i}$, where $\gamma^{\kappa}_i\in(0,1)$ and $L^{(1,\kappa)}_i$, $L^{(2,\kappa)}_i$ are positive constants. In other words,
for, $y=1,2,\ldots$,
$${\sf P}\left(Y^{(j,\kappa)}_i=y\right)=L^{(1,\kappa)}_i\left(\exp\left\{-L^{(2,\kappa)}_i
(y-1)^{\gamma^{\kappa}_i}\right\}- \exp\left\{-L^{(2,\kappa)}_iy^{\gamma^{\kappa}_i}\right\}\right),$$
or equivalently
\begin{equation}\label{example_tails}
{\sf P}\left(Y^{(j,\kappa)}_i>y\right)=L^{(1,\kappa)}_i\exp\left\{-L^{(2,\kappa)}_iy^{\gamma^{\kappa}_i}\right\},\quad y\in\mathbb{Z}_+.
\end{equation}
These formulas are particular cases of (\ref{assumption:tails_right_semiexp}) and (\ref{assumption:tails_left_semiexp}) with functions $L^{(1,\kappa)}_i(y)=L^{(1,\kappa)}_i$ and $L^{(2,\kappa)}_i(y)=L^{(2,\kappa)}_i$ for any $y\in\mathbb{Z}_+$, each $\kappa\in\{+,-\}$ and $i=1,\ldots,d$.

We verify the validity of assumption (\ref{assumptions:large_deviations_tails_right}), without loss of generality, for ${\bf x}\in\mathbb{R}^d_+$, ${\bf x}\neq{\bf 0}$. Choose parameters $L^{(1,\kappa)}_i$, $\kappa\in\{+,-\}$, for each fixed $i=1,\ldots,d$, in such a way that ${\sf E}Y^j_i=0$. Then relation (\ref{example1}) and Theorem 5.4.1 in \cite{Borovkov_Borovkov_08} imply
$${\sf P}_{\bf 0}\left({\bf S}(u)\in\Delta({\bf x};t)\right)=\sum_{k=1}^{\infty}{\sf P}\left(\Pi(u)=k\right)\prod_{i=1}^d k\left({\sf P}\left(Y^j_i\geq x_iR^{-1,+}_i(t)\right)\right)^{\left(1+\varepsilon_i(t)\right)}(1+\delta_i(k,t))$$
$$=\left(\sum_{k=1}^{\infty}k^d{\sf P}\left(\Pi(u)=k\right)\prod_{i=1}^d\left(1+\delta_i(k,t)\right)\right)\prod_{i=1}^d{\left({\sf P}\left(Y^j_i\geq x_iR^{-1,+}_i(t)\right)\right)^{\left(1+\varepsilon_i(t)\right)}},$$
where $\varepsilon_i(t)\to0$, as $t\to\infty$, and $\delta_i(k,t)\to0$, as $k,t\to\infty$. The latter equality entails the desired formula (\ref{assumptions:large_deviations_tails_right}), where $h(u)(1+\delta(u,t))=e^{-qu}\sum_{k=1}^{\infty}k^d(qu)^k/k!\prod_{i=1}^d{\left(1+\delta_i(k,t)\right)}\sim (qu)^d$, as $u,t\to\infty$, $u\leq t$.

Thus, whenever hypothesis (\ref{example_tails}) holds true for supercritical CBRW, all the conditions of Theorem~\ref{T:main_result_semiexp_lattice} are satisfied. As mentioned above, in this case we assume that $L^{(1,+)}_i+L^{(1,-)}_i=1$ and ${\sf E}Y^j_i=0$. It means that parameters $L^{(1,\kappa)}_i$, $\kappa\in\{+,-\}$, satisfy, for each $i=1,2$, the following equations system
$$\left\{
      \begin{aligned}
        &L^{(1,+)}_i+L^{(1,-)}_i=1,\\
        &L^{(1,+)}_i\sum_{y=0}^{\infty}{\exp\left\{-L^{(2,+)}_iy^{\gamma^+_i}\right\}}
-L^{(1,-)}_i\sum_{y=0}^{\infty}{\exp\left\{-L^{(2,-)}_iy^{\gamma^-_i}\right\}}=0.
      \end{aligned}
    \right.
$$
For each $i=1,2$, we have two unknown variables $L^{(1,\kappa)}_i$, $\kappa\in\{+,-\}$,
and two relations involving them in the latter system.
For example, focusing on the case $d=2$ and setting $\gamma^+_1=3/4$, $\gamma^+_2=1/2$, $\gamma^-_1=1/3$, $\gamma^-_2=1/4$, $L^{(2,+)}_1=1$, $L^{(2,+)}_2=2$, $L^{(2,-)}_1=3$ and $L^{(2,-)}_2=4$, we solve the corresponding systems with the help of packet Wolfram Mathematica and find $L^{(1,+)}_1\approx 0.382737$, $L^{(1,-)}_1\approx 0.617263$ and $L^{(1,+)}_2\approx 0.450655$, $L^{(1,-)}_2\approx 0.549345$.

Nevertheless, the limiting shape $\mathcal{P}$ of the front of the particles population described in (\ref{def_P_semiexp}) is determined exceptionally by parameters $\gamma^{\kappa}_i$, $\kappa\in\{+,-\}$, $i=1,\ldots,d$, and the Malthusian parameter $\nu$. So, to compare different forms of $\mathcal{P}$, we do no need to specify other parameters such as $L^{(1,\kappa)}_i$ and $L^{(2,\kappa)}_i$, $\kappa\in\{+,-\}$, $i=1,\ldots,d$.

\begin{center}
\begin{figure}
\includegraphics[width=17cm]{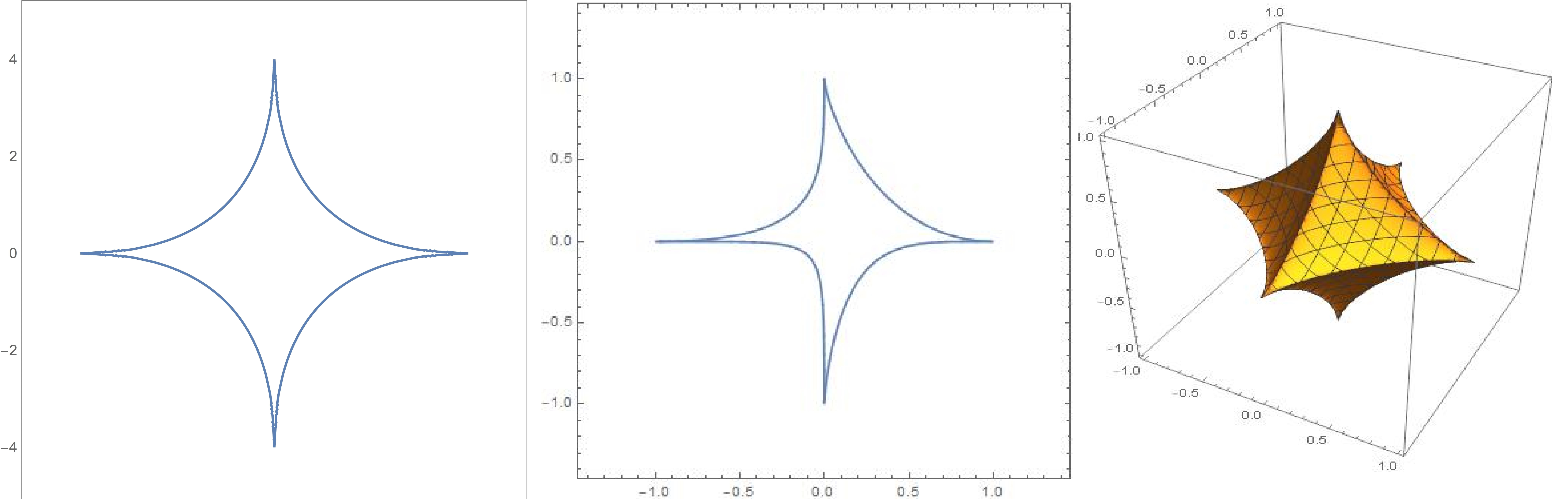}
\caption[]{The plots of $\mathcal{P}$ for Examples~1, 2 and 3.}\label{Example_123}
\end{figure}
\end{center}

\vskip0.2cm {\it Example 1.} Let $d=2$ and put $\gamma^+_1=\gamma^+_2=\gamma^-_1=\gamma^-_2=1/2$, $\nu=2$. Then the plot of the limiting shape $\mathcal{P}$ of the front is drawn on Figure~\ref{Example_123} to the left.

\vskip0.2cm {\it Example 2.} Consider now non-symmetric limiting shape $\mathcal{P}$ of the front from the example above with $d=2$, $\gamma^+_1=3/4$, $\gamma^+_2=1/2$, $\gamma^-_1=1/3$, $\gamma^-_2=1/4$, and $\nu=1$. Its plot is represented on Figure~\ref{Example_123} at the middle.

\vskip0.2cm {\it Example 3.} For $d=3$ and $\gamma^{\kappa}_i=3/4$, $\kappa\in\{+,-\}$, $i=1,2,3$, $\nu=1$, the plot of $\mathcal{P}$ is drawn on Figure~\ref{Example_123} to the right.

{\bf Remark 4}.
These Figures illustrate the fact that for CBRW on $\mathbb{Z}^d$ with semi-exponential increments the surface $\mathcal{P}$ is a boundary of a \emph{star shape set} in $\mathbb{R}^d$ with the center at ${\bf 0}$. The set is non-convex, though radially-convex. This arises from our additional assumption (\ref{assumptions:large_deviations_tails_right}) implicitly comprising
the condition of independence of the random walk jump vector coordinates. As a consequence,
to reach a distant set along the semiaxes is more probable than to reach a distant set, say,
diagonal-wise. Indeed, in the first case it is enough to perform one ``big jump'' whereas in
the second case we need to perform several ``big'' jumps along different semiaxes. However,
in the case of light tails large deviations are due to many ``small'' jumps rather than one or few ``big'' jumps, \cite{Borovkov_Borovkov_08}, p.XX.

\section{Conclusion}\label{s:conclusion}
Theorem~\ref{T:main_result_semiexp_lattice} is a counterpart of Theorem 1.1 in \cite{Carmona_Hu_14} for $d=1$ and Theorem~1 in \cite{B_SPA_18} for $d>1$ describing the population front in CBRW on $\mathbb{Z}^d$ in the case of light-tailed
jump distribution of the random walk. The novelty of our results is the following.

\begin{enumerate}
  \item The normalization of the vector ${\bf X}^v(t)$, determining the almost sure limiting shape of the front, in general is defined component-wise. It
explicitly depends both on the sign of the component
and the jump distribution of the random walk in the corresponding direction, whereas in \cite{B_SPA_18} the normalizing
factor is the same for all components of ${\bf X}^v(t)$ and just equals $t$.
\item The asymptotic behavior of the normalizing factors for the front of CBRW on $\mathbb{Z}^d$ with semi-exponential distribution tails takes an intermediate position between CBRW with light and regularly varying tails
      since the normalization of each component grows as a regularly varying function of index exceeding $\mathrm{1}$. Hence the growth is faster than linear (as for light tails). In the case of regularly varying tails
      the normalizing factor grows exponentially fast as shown in \cite{B_Arxiv_18}.
      The reason is that semi-exponential distribution tails take an intermediate position between the light and the regularly varying tails.
  \item Pictures in Section~\ref{s:examples_semiexp}
illustrate that the limiting shape $\mathcal{P}$ of the front is a boundary of a nonconvex star shape set in $\mathbb{R}^d$ which sharply contrasts to
limiting shape of the front in the case of light tails (where one has a boundary of a convex set). The nature of this effect
is explained in Remark 4.
      \end{enumerate}

We have demonstrated that for the new CBRW model it is possible to determine the limiting shape of the front for the appropriately rescaled positions of particles population. So, one can imagine that the spread of the cloud of particles  is described in time by a surface $\mathcal{P}$ of points with their coordinate components multiplied by the explicitly indicated normalizing factors (in general, different for each component). It would be interesting to study the fluctuations of
particles around this moving surface in $\mathbb{R}^d$ ($d>1$) or, equivalently, it means the analysis of the convergence rate (in a sense) of the normalized particles positions around the surface $\mathcal{P}$.

\end{document}